\documentclass[11pt]{article}
\usepackage[margin=1.5in]{geometry}
\usepackage{amssymb}
\usepackage{enumitem}
\usepackage{amsmath}
\usepackage{amsthm}
\usepackage{graphicx}
\usepackage{caption}
\usepackage{siunitx}
\usepackage{mathrsfs}
\usepackage{quiver}
\usepackage{tikz}
\usepackage{tikz-cd}
\usepackage{adjustbox}
\usepackage{hyperref}
\usepackage{bbold}
\usepackage{qworld}
\usepackage{qworld-preview}
\usepackage[affil-it]{authblk}

\DeclareMathOperator{\Hom}{Hom}

\DeclareMathOperator{\Ch}{Ch}
\DeclareMathOperator{\SBim}{SBim}

\DeclareMathOperator{\pTr}{pTr}

\newcommand{\qcdots}{\qbox[vlen=1,frame color=none,n=0,s=0,name={$\cdots$}]}
\newcommand{\qpolestick}[1][black]{\qspace \qwire[color=red] \qspace[-2] \qspider[color=black,n=3,s=3,N={1},S={1},cod arrowtype={#1},dom arrowtype={#1}] \qspace[-1]}

\newtheorem{thm}{Theorem}[section]
\newtheorem{proposition}[thm]{Proposition}

\newtheorem{lemma}[thm]{Lemma}

\theoremstyle{definition}
\newtheorem{definition}{Definition}[section]

\theoremstyle{remark}

\numberwithin{equation}{section}

\usetikzlibrary{calc, shapes, backgrounds, positioning}

\title{The partial trace functor on Soergel bimodules in type $\mathbb{B}$}
\author{Nur \c{C}ataltepe\thanks{\texttt{cataltepe.a@northeastern.edu}}}
\affil{Department of Mathematics, Northeastern University}
\date{November 19, 2026}
\begin{document}
\maketitle
\begin{abstract}
\noindent We write down some observations about the partial trace functor
between the categories of Soergel bimodules of type $\mathbb{B}_{n+1}$ and
$\mathbb{B}_n$ for future reference. In particular, we observe that, like in type $\mathbb{A}$,
the partial trace sends Soergel bimodules to a semi-infinite $q$-graded sum
of Soergel bimodules.
\end{abstract}
This serves as a reference text collating some already-established
folkloric knowledge about the
construction and
behavior of the partial trace
\[\pTr_n:\Ch^-(\SBim_{\mathbb{B}_{n+1}})\to
\Ch^-(\SBim_{\mathbb{B}_{n}}).\]
For convenience, we will be using alternate grading symbols
$q:=Q^2$ and $t:=T^2Q^{-2}$, where $T$ is the usual cohomological shift functor
and $Q$ is the usual Soergel bimodule grading shift functor. We refer
the reader to the textbook~\cite{elias_introduction_2020} for a general reference on
Soergel bimodules and the diagrammatics of the type $\mathbb{B}$ Coxeter
system.
\\ \\
We write agnostically
with respect to what realization we use for our Coxeter system, beyond the
fact that the realization of $\mathbb{B}_{n+1}$ should be one rank higher than that
of $\mathbb{B}_n$ (via the addition of an additional basis vector $x_n$, giving
an additional variable in our symmetric algebra).
\subsection*{Acknowledgements}
This paper benefited greatly from feedback and advice from the author's Ph.D.
advisor, Matthew Hogancamp.
\\ \\
\textbf{AI disclosure:} Generative AI was not used when writing this paper.
\section{The decategorified partial trace for the type $B$ Hecke algebra}
The partial trace from $\Ch^-(\SBim_{\mathbb{B}_{n+1}})$ to $\Ch^-(\SBim_{\mathbb{B}_{n}})$ is a
categorified notion of the partial trace of the Hecke algebra
$\mathbf{H}_{\mathbb{B}_{n+1}}$ when viewed as a $\mathbf{H}_{\mathbb{B}_n}$-bimodule. In this
section, we will be working with the costandard basis $\{\delta_w\}_{s\in W}$
of the Hecke algebra, which are products of the elementary costandards
$\delta_s$ for $s\in S$. When convenient, we will write $\delta_{i}$ to mean
$\delta_{s_i}$, given $s_i\in S$.
\begin{lemma}
	\label{lem:algebradecomposition}
	When viewed as a $\mathbf{H}_{\mathbb{B}_n}$-bimodule, $\mathbf{H}_{\mathbb{B}_{n+1}}$
	decomposes as
	\[\mathbf{H}_{\mathbb{B}_{n+1}}\simeq
		\mathbf{H}_{\mathbb{B}_n}\oplus \mathbf{H}_{\mathbb{B}_n}\delta_{n+1}\mathbf{H}_{\mathbb{B}_n}
		\oplus \mathbf{H}_{\mathbb{B}_n}\delta_{n+1}\dots \delta_1\dots \delta_{n+1}
	\mathbf{H}_{\mathbb{B}_{n}}.\]
Note that the last summand can also be writen omitting the
$\mathbf{H}_{\mathbb{B}_{n}}$ on the right, as $\delta_{n+1}\dots \delta_1\dots \delta_{n+1}$
commutes with all elements of $\mathbf{H}_{\mathbb{B}_{n}}$ and its corresponding
element in $W$ only has $s_{n+1}$ in either its left or right descent set.
\end{lemma}
\begin{proof}
	For $n=1$, the Coxeter group of type $\mathbb{B}_2$ consists of
\[\{1,s_1,s_2,s_1s_2,s_2s_1,s_1s_2s_1,s_2s_1s_2,s_1s_2s_1s_2\}\]
where the sets of corresponding costandard basis elements
\[\{1,\delta_1\},\ \{\delta_2,\delta_2\delta_1,
\delta_1\delta_2,\delta_1\delta_2\delta_1\},\ \{\delta_2\delta_1\delta_2,
\delta_2\delta_1\delta_2\delta_1\}\]
obviously span disjoint sub-$\mathbf{H}_{\mathbb{B}_1}$-bimodules.
\\ \\
We now assume $n \geq 2$.
	The longest element of the Coxeter group of
	type $\mathbb{B}_{n+1}$ is of the form
	\[\prod_{i=1}^{n+1}s_is_{i-1}\dots s_1 \dots s_{i-1}s_i \]
	and since every element of the Coxeter group appears as a sub-expression of
	the longest element, this means that the elements of the Coxeter
	group of type $\mathbb{B}_{n+1}$ divide into three types:
	\begin{itemize}
	\item Those which have a reduced expression in which
		$s_{n+1}$ does not appear. These are elements of the
		parabolic subgroup of type $\mathbb{B}_n$.
	\item Those which have a reduced expression in which $s_{n+1}$ appears once.
	\item Those which only have reduced expressions in which $s_{n+1}$ appears
		twice or more (and, by virtue of being realizable as
		sub-expressions of the longest element, have a reduced expression
		in which $s_{n+1}$ appears exactly twice).
	\end{itemize}
	We now show that the
	costandard basis elements corresponding to the first two of these types
	give disjoint
	sub-$\mathbf{H}_{\mathbb{B}_n}$-bimodules of $\mathbf{H}_{\mathbb{B}_{n+1}}$.
	If we take $w$ to be an element of
	the Coxeter group, this follows from the following arguments:
	\begin{itemize}
		\item If $w$ has a reduced expression in which $s_{n+1}$ does
			not appear, then as it is in a parabolic subgroup,
			$ws_i$ for $i < n+1$ will also have a reduced expression
			without $s_{n+1}$, and so $\delta_w\delta_{i}\in
			\mathbf{H}_{\mathbb{B}_n}$.
		\item If $w$ has a reduced expression in which $s_{n+1}$
			appears once, write a
			reduced expression $w=vs_{n+1}v'$ such that
			$v,v'$ do not have $s_{n+1}$ appear in them. Take some
			$s_i$ for $i < n+1$. If $vs_{n+1}v's_i$ is a reduced
			expression fo $ws_i$, then we are done (since
			$i < n+1$, this still has one instance of
			$s_{n+1}$). Otherwise,
			the exchange condition dictates that there is some
			simple reflection appearing in $vs_{n+1}v's_i$ we can
			delete to produce a reduced expression for $ws_i$.
			We can appeal to the fact that
			the strand diagram for $ws_i$ is of the form
			\[\q{
			\qwire
			\qcdots \qcdots \qcdots \qcdots \qwire[color=red]\n
				\qsym  \qcdots \qcdots \qcdots
				\qwire[color=red] \n
				\qwire \qcdots \qcdots \qcdots \qcdots
			\qwire[color=red] \n
		\qwire \qcdots \qsym \qcdots \qwire[color=red]}
	\ \mathrm{or}\ 
	\q{\qwire
			\qcdots \qcdots \qcdots \qwire[color=red]\n
				\qsym  \qcdots \qcdots
				\qwire[color=red] \n
				\qwire \qcdots \qcdots \qcdots
			\qwire[color=red] \n
		\qwire \qcdots \qcdots \qpolestick}
	\]
	and thus $s_i$ cannot share both of its strands with $s_{n+1}$, meaning
	$s_{n+1}$ cannot be the element to be deleted. Thus,
	$\delta_w\delta_{i}\in \mathbf{H}_{\mathbb{B}_n}\delta_{n+1}\mathbf{H}_{\mathbb{B}_n}$.
	\end{itemize}
	We now prove that any element $w$
	of the Coxeter group with a reduced expression
	$s_{n+1}vs_{n+1}$, where $v$ is in the parabolic subgroup of type $\mathbb{B}_n$,
	must either be expressible in the form
	$s_{n+1}s_n\dots s_1 \dots s_n s_{n+1}v'$ where
	$v'$ is in the parabolic subgroup, or have a reduced expression in which
	$s_{n+1}$ appears one or zero times.
	\\ \\
	To start with, $w$ has a string diagram of the form
\[\q{		\qsym[R color=blue] \qwire \qcdots \qwire[color=red] \n
	\qwire \qcdots \qcdots \qcdots \qcdots \n
	\qsym[L color=purple] \qwire \qcdots \qwire[color=red]

}\]
If the blue and purple strands connect while hitting the red pole, we are done, since
this creates an instance of $s_{n+1}s_n\dots s_{1}\dots s_ns_{n+1}$ which can be
moved past the rest of the expression. If they connect without hitting the red pole,
then this means the two instances of $s_{n+1}$ cancel each other out, and we
are also done.
\\ \\
Now, suppose the blue and purple strands do not connect, in which they case they
must cross. If neither strand hits the red pole, we can apply a Reidemeister III
move to obtain
\[\q{\qsym[R color=blue]  \qcdots \qwire[color=purple]\qcdots   \qwire[color=red] \n
	\qwire \qcdots \qcdots \qwire[color=purple] \qcdots \qwire[color=red] \n
		\qwire \qcdots \qsym[R color= blue,
		L color=purple] \qcdots  \qwire[color=red]\n
		\qwire \qcdots \qcdots \qwire[color=blue] \qcdots \qwire[color=red] \n
		\qsym[L color=purple] \qcdots \qwire[color=blue]
		\qcdots  \qwire[color=red]

	}=
\q{
	\qwire[color=blue] \qsym[L color=purple,hlen=3] \qcdots \qwire[color=red]\n
	\qwire[color=blue] \qcdots \qcdots \qwire \qcdots \qwire[color=red] \n
	\qsym[R color=blue, L color=purple] \qcdots \qwire \qcdots \qwire[color=red]\n
	\qwire[color=purple] \qcdots \qcdots \qwire \qcdots \qwire[color=red] \n
	\qwire[color=purple] \qsym[R color=blue,hlen=3] \qcdots \qwire[color=red]\n
}
\]
which means $w$ has a reduced expression with only one occurrence of
$s_{n+1}$.
\\\\
If either or both of the strands hits the red pole, that is, if
$w$ has a strand diagram of the form
\[
	\q{
		\qsym[R color=blue] \qcdots \qpolestick[purple] \n
		\qwire \qcdots \qcdots  \qwire[color=purple] \qwire[color=red] \n
\qwire \qcdots \qsym[R color=blue, L color=purple]\qwire[color=red] \n
\qwire \qcdots \qcdots \qwire[color=blue] \qwire[color=red] \n
\qsym[L color=purple] \qcdots \qwire[color=blue] \qwire[color=red]}\ \mathrm{or}\
\q{
		\qsym[R color=blue] \qcdots \qpolestick[purple] \n
		\qwire \qcdots \qcdots  \qwire[color=purple] \qwire[color=red] \n
\qwire \qcdots \qsym[R color=blue, L color=purple]\qwire[color=red] \n
\qwire \qcdots \qcdots \qwire[color=blue] \qwire[color=red] \n
\qsym[L color=purple] \qcdots \qpolestick[blue]}
\]
then we can use the fact that $n \geq 2$ to move the instances of $s_1$ past
$s_{n+1}$ to apply the same argument as before.
\\ \\
Now suppose $w$ to be an element of the Coxeter group of type $\mathbb{B}_{n+1}$ which
only has reduced expressions in which $s_{n+1}$ appears twice or more, and in
particular has a reduced expression in which $s_{n+1}$ appears exactly twice.
Since the longest element in
type $\mathbb{B}_{n+1}$ can be written such that $s_{n+1}$ is the first simple reflection
to appear in it, we can write a reduced expression $w=s_{n+1}vs_{n+1}v'$, where
$v,v'$ are in the parabolic subgroup. We have that $s_{n+1}vs_{n+1}$ is a
contiguous sub-expression of a reduced expression and itself reduced, so
if $s_{n+1}vs_{n+1}$ has a reduced expression in which $s_{n+1}$ only appears
once or fewer times, then so does $w$. Thus, $s_{n+1}vs_{n+1}$ must only
have reduced expressions in which $s_{n+1}$ appears twice or more, and must
thus have a reduced expression obtainable by braid moves as
$s_{n+1}\dots s_1\dots s_{n+1}v$, as must $w$. Now, $w$ has a strand
diagram of the form
\[\q{\qsym[R color = blue] \qcdots \qcdots \qcdots \qwire[color=red]\n
\qwire \qcdots \qcdots \qcdots \qpolestick[blue] \n
\qsym[L color = blue] \qcdots \qcdots \qcdots \qwire[color=red]\n
\qwire[color=blue] \qcdots \qcdots \qcdots \qcdots \qwire[color=red]}\]
Given some $s_i$ where $i < n+1$, suppose $s_{n+1}\dots s_1\dots s_{n+1}vs_i$
is a reduced expression for $ws_i$. Then, since braid relations involving
an element in the parabolic subgroup can only increase the number of $s_{n+1}$s,
$\delta_w\delta_{i}$ is still in $\delta_{n+1}\dots \delta_1\dots
\delta_{n+1}\mathbf{H}_{\mathbb{B}_n}$.
\\ \\
Now, suppose $s_{n+1}\dots s_1\dots s_{n+1}vs_i$ is not a reduced expression.
Then, there is an element of
$s_{n+1}\dots s_1\dots s_{n+1}v$ we can delete to obtain a reduced expression
for $ws_i$. However, as we can see from the strand diagram for $w$, no simple
reflection in the parabolic subgroup can connect to the blue strand and so no simple
reflection appearing in $s_{n+1}\dots s_1\dots s_{n+1}$ can be deleted. Therefore,
$\delta_w\delta_i\in \delta_{n+1}\dots \delta_1\dots \delta_{n+1}\mathbf{H}_{\mathbb{B}_n}$.
\\ \\
As $\mathbf{H}_{\mathbb{B}_n}$, $\mathbf{H}_{\mathbb{B}_n}\delta_{n+1}\mathbf{H}_{\mathbb{B}_n}$ and
$\delta_{n+1}\dots \delta_1\dots \delta_1 \mathbf{H}_{\mathbb{B}_n}$ are disjoint
sub-bimodules which together span $\mathbf{H}_{\mathbb{B}_{n+1}}$, the
direct sum decomposition follows.
\end{proof}
\noindent
We are now ready to define the decategorified partial trace.
\begin{definition}
Given the decomposition
\[\mathbf{H}_{\mathbb{B}_{n+1}}=
	\mathbf{H}_{\mathbb{B}_n} \oplus \mathbf{H}_{\mathbb{B}_n}\delta_{n+1}\mathbf{H}_{\mathbb{B}_n}
\oplus \delta_{n+1} \dots \delta_{1}\dots \delta_{n+1}\mathbf{H}_{\mathbb{B}_n}\]
the \textbf{partial trace} is the $\mathbf{H}_{\mathbb{B}_n}$-bimodule homomorphism
\[\pTr_n:\mathbf{H}_{\mathbb{B}_{n+1}}\to \mathbf{H}_{\mathbb{B}_n}\]
which is the identity on $\mathbf{H}_{\mathbb{B}_n}$ and zero on the other
two direct summands.
\end{definition}
\section{The (right) partial trace $\pTr_{n}$ on costandard objects $\nabla_i$}
We now categorify the partial trace as a $dg$-functor between the categories
of chain complexes of Soergel bimodules. We first recall the costandard objects
in these categories.
\begin{definition}
Given a Coxeter system $(W,S)$ and a simple reflection $s\in S$, the
\textbf{elementary costandard object} $\nabla_s\in
\Ch^-(\SBim_{(W,S)})$ is given by
\[
\nabla_s:=
t^{-1/2}\mathbb{1} \xrightarrow{\q[scale=0.4]{\qunit[color=black]}} B_s\]
where $B_{s}$ is the elementary Bott-Samelson bimodule corresponding to
$s$. Elementary costandard objects satisfy categorified versions of the braid
relations and the quadratic relation in the Hecke algebra up to homotopy equivalence.
The monoidal identity $\mathbb{1}$ is considered the costandard object $\nabla_1$
corresponding to the identity in $W$.
\\ \\
Given some $w\in W$ and a reduced expression $s_1\dots s_k$ for $w$, we define
the \textbf{costandard object} $\nabla_w$ corresponding to $w$ to be
\[\nabla_w:=\nabla_{s_1}\dots \nabla_{s_k}\]
which is well defined up to homotopy equivalence. 
\end{definition}
\noindent
Notice that given any elementary Bott-Samelson $B_s$, we can write a homotopy
equivalence
\[B_s\simeq \nabla_s \to q^{1/2}\mathbb{1}\]
where the twist is simply an identity between the occurrences of
$\mathbb{1}$ in either complex. Therefore, $\Ch^-(\SBim_{\mathbb{B}_{n}})$ is generated
by the collection of costandard objects
$\{\nabla_w,w\in W\}$ under grading shifts, direct sums, twists,
homotopy equivalences, and bounded-above cochain complexes.
\subsection{Categorifying the partial trace}
There is a natural ``inclusion'' functor $I_n:\SBim_{\mathbb{B}_n}\to \SBim_{\mathbb{B}_{n+1}}$, which
sends Soergel bimodules $B_w$ to the corresponding Soergel bimodules in
$\SBim_{\mathbb{B}_{n+1}}$, and sends the morphism spaces in $\SBim_{\mathbb{B}_n}$ to those
between the bimodules with the same labels in $\SBim_{\mathbb{B}_{n+1}}$ by tensoring with
$\mathbb{k}[x_{n+1}]$. When extended
to the categories of bounded above cochain
complexes, the inclusion functor sends costandards to
costandards as well.
\\ \\
We are motivated to categorify the partial trace so as to have an adjoint to the
inclusion functor. In this text, we define a partial trace on costandard
objects that acts as a right adjoint to $I_n$, but near-identical arguments
with standard objects (positive Rouquier complexes) will yield a left adjoint.
\\ \\
We will define our adjoint initially on costandard objects as an additive functor,
which will extend as a $dg$-functor to all of $\Ch^-(\SBim_{\mathbb{B}_{n+1}})$. We use the
following property of costandard objects, which
easily follows from Rouquier's formula (given at the start of Section 2.5
of~\cite{libedinsky_standard_2014}) as an exercise,
to easily categorify the partial
trace.
\begin{lemma}[\cite{libedinsky_standard_2014}]
	\label{lem:costandardhoms}
	Let $x,y$ be elements of the Coxeter group of type $\mathbb{B}_{n+1}$ such that
	$x \ngeq y$ in the Bruhat order.
	\\\\
	Then,
	given $\nabla_{x},\nabla_y\in \Ch^-(\SBim_{\mathbb{B}_{n+1}})$,
	$\Hom(\nabla_x,\nabla_y)\simeq 0$.
\end{lemma}
We can now define the partial trace functor.
\begin{definition}
	The \textbf{(right) partial trace} $\pTr_n:\Ch^-(\SBim_{\mathbb{B}_{n+1}})\to
	\Ch^-(\SBim_{\mathbb{B}_{n}})$ is the functor defined such that
\begin{itemize}
	\item $\pTr_n(\nabla_x)=\nabla_x\otimes \mathbb{k}[x_{n+1}]\in
		\Ch^-(\SBim_{\mathbb{B}_n})$ for $x$ in the parabolic subgroup of
		type $\mathbb{B}_n$.
	\item $\pTr_n(\nabla_x)=0$ otherwise.
\end{itemize}
\end{definition}
\noindent
That this is a functor is obvious, and that it provides a right adjoint to
$I_n$ is immediate from Lemma~\ref{lem:costandardhoms}, as $x\ngeq y$ if
$x$ is in the parabolic subgroup of type $\mathbb{B}_n$ but $y$ is not. Stated more
formally,
\begin{proposition}
	\label{prop:adjoint}
	The right partial trace $\pTr_n$ is a right adjoint to $I_n$, that is,
	given $X\in \Ch^-(\SBim_{\mathbb{B}_n})$ and
	$Y\in \Ch^-(\SBim_{\mathbb{B}_{n+1}})$,
	there is a natural homotopy equivalence
	\[\Hom_{\Ch^-(\SBim_{\mathbb{B}_{n+1}})}
	(I_n(X),Y)\simeq \Hom_{\Ch^-(\SBim_{\mathbb{B}_n})}(X,\pTr_n(Y)).\]
\end{proposition}
\subsection{Monoidality}
The decategorified partial trace is a $\mathbf{H}_{\mathbb{B}_n}$-bimodule
homomorphism. The category $\Ch^-(\SBim_{\mathbb{B}_n})$ monoidally acts on
$\Ch^-(\SBim_{\mathbb{B}_{n+1}})$ via $I_n$ and tensoring, and it can easily be
seen that $\pTr_n$ intertwines with this action.
\begin{proposition}
	\label{prop:monoidality}
	The right partial trace $\pTr_n$ is a functor of categorical modules
	of $\Ch^-(\SBim_{\mathbb{B}_n})$, that is, given $X\in \Ch^-(\SBim_{\mathbb{B}_n})$
	and $Y\in \Ch^-(\SBim_{\mathbb{B}_{n+1}})$, there is a natural homotopy equivalence
	\[\pTr_n(I_n(X)Y)\simeq X\pTr_n(Y).\]
\end{proposition}
\section{The partial trace on Soergel bimodules}
In the case of type $\mathbb{A}$, it has proven a convenient fact that the partial trace
restricts to a functor between Soergel bimodules (up to a factor of
$\frac{1}{1-q}$). We show an analogous result for type $\mathbb{B}$.
\begin{proposition}
	The right partial trace $\pTr_n$ sends Soergel bimodules in type $\mathbb{B}_{n+1}$
	to Soergel bimodules in type
	$\mathbb{B}_n$ with a factor of $\frac{1}{1-q}$ (tensoring with
	$\mathbb{k}[x_{n+1}]$).
\end{proposition}
\begin{proof}
	Due to the decomposition proven in Lemma~\ref{lem:algebradecomposition},
	we have $3$ cases for indecomposable Soergel bimodules:
	\begin{itemize}
		\item $B_x$ for $x$ in the parabolic subgroup of type $\mathbb{B}_n$. Since every
		costandard in the costandard resolution is simply sent to itself
		times $\frac{1}{1-q}$, so is $B_x$.
	\item $B_{xs_{n+1}y}$, with $x,y$ in the parabolic subgroup. As
		$\pTr_n$ is an additive functor, $\pTr_n(B_{xs_{n+1}y})$ is
		a direct summand of
		$\pTr_n(B_xB_{n+1}B_y)$, which is homotopy equivalent to
		\[\pTr(B_x(\nabla_{n+1}\to q^{1/2}\mathbb{1})B_y)
		=\pTr(B_x\nabla_{n+1}B_y) \to q^{1/2}\pTr(B_xB_y).\]
		As the minimal
		complex for the
		standard resolution of $B_x\nabla_{n+1}B_y$ will have
		a factor of $\nabla_{n+1}$ in every term, it will be sent to
		zero, so the partial trace will be that of $q^{1/2}B_xB_y$.
		Therefore, the partial trace of our Soergel bimodule will be
		a direct summand of $\frac{1}{1-q}$ 
	\item $B_{s_{n+1}\dots s_1\dots s_{n+1}x}$, for $x$ in the parabolic
		subgroup.
	\end{itemize}
	For the last case, we proceed by induction on $n$. For $n=1$,
	observe that the standard resolution for $B_{212}$ is
	\[\begin{tikzcd}[ampersand replacement=\&,cramped]
	\& {q^{1/2}\nabla_1\nabla_2} \& {q\nabla_2} \& \\
	{\nabla_{2}\nabla_1\nabla_2} \&\&\& {q^{3/2}\mathbb{1}} \\
	\& {q^{1/2}\nabla_2\nabla_1} \& {q\nabla_1}
	\arrow[from=1-2, to=1-3]
	\arrow[from=1-2, to=3-3]
	\arrow[from=1-3, to=2-4]
	\arrow[from=2-1, to=1-2]
	\arrow[from=2-1, to=3-2]
	\arrow[from=3-2, to=1-3]
	\arrow[from=3-2, to=3-3]
	\arrow[from=3-3, to=2-4]
\end{tikzcd}\]
	which will become $\frac{q}{1-q}B_1$ after applying the partial trace.
	We now recall the categorified quadratic relation for
	costandards
	\[\nabla_i^2\simeq t^{-1/2}\nabla_i\to q^{-1/2}\nabla_i\to \mathbb{1}\]
	and take $B_{32123}=B_3B_{212}B_3$, which we can partially write
	in terms of costandards as
	\[\begin{tikzcd}[ampersand replacement=\&,cramped]
	\& {\nabla_3B_{212}\nabla_3} \& \\
	{q^{1/2}\nabla_3B_{212}} \&\& {q^{1/2}B_{212}\nabla_3} \\
	\& {qB_{212}}
	\arrow[from=1-2, to=2-1]
	\arrow[from=1-2, to=2-3]
	\arrow[from=2-1, to=3-2]
	\arrow[from=2-3, to=3-2]
\end{tikzcd}\]
	The costandard resolution of the middle two terms will consist entirely
	of terms containing a $\nabla_3$, so the partial trace will send these
	to zero, and send $qB_{212}$ to $\frac{q}{1-q}B_{212}$.
	Now, for
	the term $\nabla_3B_{212}\nabla_3$, we can observe that any
	terms in the minimal costandard resolution of $B_{212}$ with a
	$\nabla_2$ will not allow the two instances of $\nabla_3$ to create
	a quadratic relation, and the terms which do allow a quadratic
	relation to occur will precisely be those lacking a factor of
	$\nabla_2$, so the partial trace of this term will be
	$I_1(\pTr_{1}(B_{212}))=\frac{q}{1-q}B_{1}$. Thus,
	\[\pTr_2(B_3B_{212}B_3)=\frac{q}{1-q}(B_{1}\oplus B_{212}).\]
	Now, assume for $n \geq 3$ that
	\begin{align*}\pTr_{n-1}(B_{s_{n}\dots s_1\dots s_{n}})&=
		\pTr_{n-1}(B_{n}\dots B_{212}\dots B_{n})\\ &=
\frac{q}{1-q}(B_{1}\oplus B_{212}\oplus \dots \oplus B_{n-1}\dots B_{212}\dots
B_{n-1})\end{align*}
and take
\[B_{n+1}B_n\dots B_{212}\dots B_nB_{n+1}\]
which we can write as
\[
\begin{tikzcd}[ampersand replacement=\&,cramped, column sep=tiny]
	\& {\nabla_{n+1}B_{n}\dots B_{212}\dots B_{n}\nabla_{n+1}} \& \\
	{q^{1/2}\nabla_{n+1}B_{n}\dots B_{212}\dots B_{n}} \&\& {q^{1/2}B_{n}\dots B_{212}\dots B_{n}\nabla_{n+1}} \\
	\& {qB_{n}\dots B_{212}\dots B_{n}}
	\arrow[from=1-2, to=2-1]
	\arrow[from=1-2, to=2-3]
	\arrow[from=2-1, to=3-2]
	\arrow[from=2-3, to=3-2]
\end{tikzcd}
\]
By arguments previously given, the middle two terms will vanish under the partial
trace, and the only terms in the top term that willl survive will be the ones
in which $\nabla_{n+1}$ can cause a quadratic relation to occur (i.e., those
with no factors of $\nabla_n$). Thus, the partial trace is
\begin{align*}\pTr(B_{n+1}\dots B_{212}\dots B_{n+1})&=
	\frac{q}{1-q}B_{n}\dots B_{212}\dots B_{n} \oplus I_{n-1}
\pTr_{n-1}(B_{n}\dots B_{212}\dots B_n)\\
	&= \frac{q}{1-q}(B_1\oplus B_{212}\oplus \dots
\oplus B_{n}\dots B_{212}\dots B_n)\end{align*}
proving the statement.
\\ \\
Since $\pTr_n(B_{s_{n+1}\dots s_1\dots s_{n+1}x})$ is a direct summand
of $\pTr_{n}(B_{s_{n+1}\dots s_1\dots s_{n+1}}B_x)$, we can use the fact that
$\pTr_n$ is a functor of categorical modules to move the $B_x$ outside the
partial trace and see that the partial trace of this element will be a direct
summand of a tensor product of Soergel bimodules and therefore a Soergel
bimodule itself. Thus, $\pTr_n$ sends Soergel bimodules to Soergel bimodules with
a factor of $\frac{1}{1-q}$.
\end{proof}
\bibliographystyle{abbrv}
\bibliography{main}
\end{document}